\documentclass[reqno]{amsart}
\usepackage{amssymb,hyperref}

\theoremstyle{plain}
\newtheorem{thm}{Theorem}
\newtheorem{lem}{Lemma}
\newtheorem{cor}{Corollary}

\theoremstyle{definition}
\newtheorem{rem}{Remark}

\renewcommand{\Re}{\mathrm{Re}}

\title
[Some sufficient problems for certain univalent functions]
{Some sufficient problems \\
for certain univalent functions}

\author{Hitoshi Shiraishi}
\address{Hitoshi Shiraishi \newline
Department of Mathematics \newline
Kinki University \newline
Higashi-Osaka, Osaka 577-8502, Japan}
\email{0733310104x@kindai.ac.jp}

\author{Shigeyoshi Owa}
\address{Shigeyoshi Owa \newline
Department of Mathematics \newline
Kinki University \newline
Higashi-Osaka, Osaka 577-8502, Japan}
\email{owa@math.kindai.ac.jp}


\date{}

\begin{document}

\begin{abstract}
For analytic functions $f(z)$ in the open unit disk $\mathbb{U}$ with $f(0)=f'(0)-1=0$,
R. Singh and S. Singh (Coll. Math. {\bf 47}(1982), 309-314) have considered some sufficient problems for $f(z)$ to be univalent in $\mathbb{U}$.
The object of the present paper is to discuss some sufficient problems for $f(z)$ to be some classes of analytic functions in $\mathbb{U}$.
\end{abstract}

\begin{flushleft}
This paper was published in the journal: \\
Far East J. Math. Sci. (FJMS) {\bf 30} (2008), No. 1, 147--155.
\end{flushleft}
\hrule

\

\

\maketitle

\section{Introduction}

\

Let $\mathcal{A}$ denote the class of functions $f(z)$ that are analytic in the open unit disk $\mathbb{U}=\{z \in \mathbb{C}:|z|<1\}$,
so that $f(0)=f'(0)-1=0$.

We denote by $\mathcal{S}$ the subclass of $\mathcal{A}$ consisting of univalent functions $f(z)$ in $\mathbb{U}$.

\

Let $\mathcal{C}(\alpha)$ denote
$$
\mathcal{C}({\alpha})
= \{f(z)\in\mathcal{A}:|f'(z)-1| < 1-\alpha,\ 0\leqq\alpha<1\}
$$
and $\mathcal{C}=\mathcal{C}(0)$.

Also, let $\mathcal{S^{*}}(\alpha)$ be defined by
$$
\mathcal{S^{*}}(\alpha)
= \left\{f(z)\in\mathcal{A}:\Re\left(\frac{zf'(z)}{f(z)}\right) > \alpha,\ 0\leqq\alpha<1\right\}
$$
and $\mathcal{S^{*}}=\mathcal{S^{*}}(0)$.

Further, let $\mathcal{STS}(\mu)$ denote
$$
\mathcal{STS}(\mu)
= \left\{f(z)\in\mathcal{A}:\Re\left(\frac{zf'(z)}{f(z)}\right)^{\frac{1}{\mu}}>0,\ 0<\mu\leqq1\right\}
$$
and $\mathcal{STS}=\mathcal{STS}(1)$.

\

The basic tool in proving our results is the following lemma due to Jack \cite{m1ref1} (also, due to Miller and Mocanu \cite{m1ref2}).

\

\begin{lem} \label{jack} \quad
Let $w(z)$ be analytic in the open unit disk $\mathbb{U}$ with $w(0)=0$. Then if $\left|w(z)\right|$ attains its maximum value on the circle $\left|z\right|=r$ at a point $z_{0}\in\mathbb{U}$, then we have $z_{0}w'(z_{0})=kw(z_{0})$, where $k\geqq1$ is a real number.
\end{lem}

\

\section{Conditions for the class \texorpdfstring{$\mathcal{C}$}{C}}

\

Applying Lemma \ref{jack},
we drive the following result for the class $\mathcal{C}$.

\

\begin{thm} \label{m1thm1} \quad
If $f(z)\in\mathcal{A}$ satisfies
$$
\left|f'(z)-1\right|^{\beta}\left|\delta+\frac{zf''(z)}{f'(z)}\right|^{\gamma}<\left(\frac{1+2\delta}{2}\right)^{\gamma}
\qquad (z\in\mathbb{U})
$$
for some real $\beta$, $\gamma \geqq 0$ and $\delta > -\dfrac{1}{2}$,
then $f(z) \in \mathcal{C}$.
\end{thm}

\

\begin{proof}\quad
Let us define $w(z)$ by
\begin{equation}
w(z)=f'(z)-1 \qquad (z\in\mathbb{U}).
\label{m1thm1eq1}
\end{equation}

Then, clearly, $w(0)=0$ and $w(z)$ is analytic in $\mathbb{U}$.
Differentiating both sides in (\ref{m1thm1eq1}),
we obtain
$$
\frac{zf''(z)}{f'(z)}=\frac{zw'(z)}{1+w(z)},
$$
and therefore,
$$
\left|f'(z)-1\right|^{\beta} \left| \delta + \frac{zf''(z)}{f'(z)}\right|^{\gamma}=\left| w(z) \right| ^{\beta} \left| \delta + \frac{zw'(z)}{1+w(z)} \right| ^{\gamma} < \left(\frac{1+2\delta}{2}\right)^{\gamma}
\qquad (z\in\mathbb{U}).
$$

If there exists a point $z_{0} \in \mathbb{U}$ such that
$$
\max_{\left| z \right| \leqq \left| z_{0} \right|} \left| w(z) \right|
= \left| w(z_{0}) \right|
= 1,
$$
then Lemma \ref{jack} gives us that $w(z_{0})=e^{i \theta}$ and $z_{0}w'(z_{0})=kw(z_{0})$ $(k \geqq 1)$.

Thus we have
\begin{align*}
\left|f'(z_{0})-1\right|^{\beta} \left|\delta+\frac{z_{0}f''(z_{0})}{f'(z_{0})}\right|^{\gamma}
&= \left|w(z_{0})\right| ^{\beta} \left|\delta+\frac{z_{0}w'(z_{0})}{1+w(z_{0})}\right|^{\gamma}\\
&= \left|\delta+\frac{kw(z_{0})}{1+w(z_{0})}\right|^{\gamma}\\
&= \left|\delta+\frac{ke^{i\frac{\theta}{2}}}{e^{i\frac{\theta}{2}}+e^{-i\frac{\theta}{2}}}\right|^{\gamma}\\
&= \left(\frac{1}{2}\right)^{\gamma} \left(\left(k+2\delta\right)^{2}+k^{2}\tan^2\frac{\theta}{2}\right)^{\frac{\gamma}{2}}\\
&\geqq \left(\frac{1+2\delta}{2}\right)^{\gamma}.
\end{align*}

This contradicts our condition in the theorem.
Therefore,
there is no $z_{0} \in \mathbb{U}$ such that $\left|w(z_{0})\right| = 1$.
This means that $\left|w(z)\right| < 1$ for all $z \in \mathbb{U}$.
It follows that $\left|f'(z)-1\right| < 1\ (z \in \mathbb{U}$) so that,
$f(z)\in \mathcal{C}$.
\end{proof}

\

Letting $\beta=1-\lambda$, $\gamma=\lambda$ and $\delta=1$ in Theorem \ref{m1thm1},
we have the following corollary by Singh and Singh \cite{m1ref3}.

\

\begin{cor} \label{m1cor1} \quad
If $f(z)\in\mathcal{A}$ satisfies
$$
\left|f'(z)-1\right|^{1-\lambda}\left|1+\frac{zf''(z)}{f'(z)}\right|^{\lambda}
<\left(\frac{3}{2}\right)^{\lambda}
\qquad (z\in\mathbb{U})
$$
for some real $\lambda \geqq 0$,
then $f(z) \in \mathcal{C}$.
\end{cor}

\

Making $\delta=0$ in Theorem \ref{m1thm1},
we see

\

\begin{cor} \label{m1cor2} \quad
If $f(z)\in\mathcal{A}$ satisfies
$$
\left|f'(z)-1\right|^{\beta}\left|\frac{zf''(z)}{f'(z)}\right|^{\gamma}<\left(\frac{1}{2}\right)^{\gamma}
\qquad (z\in\mathbb{U})
$$
for some real $\beta$ and $\gamma \geqq 0$,
then $f(z) \in \mathcal{C}$.
\end{cor}

\

\begin{rem} \label{m1re1} \quad
If we take $\gamma=0$ in Corollary \ref{m1cor2},
then we have that, for some real $\beta$,
$$
\left|f'(z)-1\right|^{\beta}<1
\qquad (z\in\mathbb{U})
$$
implies
$$
\left|f'(z)-1\right|<1
\qquad (z\in\mathbb{U}).
$$
\end{rem}

\

\section{Conditions for the class \texorpdfstring{$\mathcal{S^{*}}(\alpha)$}{S*(alpha)}}

\

Next,
we derive the following result for the class $\mathcal{S^{*}}(\alpha)$.

\

\begin{thm} \label{m1thm2} \quad
If $f(z)\in\mathcal{A}$ satisfies
\begin{equation}
\left|\frac{zf'(z)}{f(z)}-1\right|^{\beta} \left|z\left(\frac{zf'(z)}{f(z)}\right)'\right|^{\gamma}
<\left(\dfrac{1}{2}\right)^{\gamma}
\qquad (z\in\mathbb{U}) \label{m1thm2eq1}
\end{equation}
or
\begin{equation}
\left|\frac{zf'(z)}{f(z)}+1\right|^{\beta} \left|z\left(\frac{zf'(z)}{f(z)}\right)'\right|^{\gamma}
< \left(\dfrac{1}{2}\right)^{\gamma}
\qquad (z\in\mathbb{U}) \label{m1thm2eq2}
\end{equation}
for some real $\beta$, $\gamma$ with $\beta + 2\gamma \geqq 0$,
then $f(z) \in \mathcal{S^{*}}$.
\end{thm}

\

\begin{proof}\quad
Define $w(z)$ in $\mathbb{U}$ by
\begin{equation}
G(z)=\frac{zf'(z)}{f(z)}=\frac{1+w(z)}{1-w(z)}
\qquad (w(z) \neq 1). \label{m1thm2eq3}
\end{equation}

Evidently,
$w(0)=0$ and $w(z)$ is analytic in $\mathbb{U}$.
Differentiating (\ref{m1thm2eq3}) logarithmically and simplifying,
we obtain
$$
\left(\frac{zf'(z)}{f(z)}\right)'=\frac{2w'(z)}{(1-w(z))^{2}}
$$
and, hence
$$
\left|\frac{zf'(z)}{f(z)}-1\right|^{\beta} \left|z\left(\frac{zf'(z)}{f(z)}\right)'\right|^{\gamma} = \left|\frac{2w(z)}{1-w(z)}\right|^{\beta} \left|\frac{2zw'(z)}{(1-w(z))^{2}}\right|^{\gamma} < \left(\dfrac{1}{2}\right)^{\gamma}
\qquad (z\in\mathbb{U}).
$$

If there exists a point $z_{0} \in \mathbb{U}$ such that
$$
\max_{\left| z \right| \leqq \left| z_{0} \right|} \left| w(z) \right|
= \left| w(z_{0}) \right|
= 1,
$$
then Lemma \ref{jack} gives us that $w(z_{0})=e^{i \theta}$ and $z_{0}w'(z_{0})=kw(z_{0})$ $(k \geqq 1)$.

Thus we have
\begin{align*}
\left|\frac{z_{0}f'(z_{0})}{f(z_{0})}-1\right|^{\beta} \left|z_{0}\left(\frac{z_{0}f'(z_{0})}{f(z_{0})}\right)'\right|^{\gamma}
&= \left|\frac{2w(z_{0})}{1-w(z_{0})}\right|^{\beta} \left|\frac{2z_{0}w'(z_{0})}{(1-w(z_{0}))^{2}}\right|^{\gamma}\\
&= \frac{2^{\beta+\gamma}k^{\gamma}}{\left|1-w(z_{0})\right|^{\beta+2\gamma}}\\
&\geqq \left(\frac{k}{2}\right)^{\gamma} \geqq \left(\frac{1}{2}\right)^{\gamma}.
\end{align*}

This contradicts the condition (\ref{m1thm2eq1}) in the theorem.
Therefore, there is no $z_{0} \in \mathbb{U}$ such that $|w(z_{0})|=1$.
This means that $\left|w(z)\right| < 1$ for all $z \in \mathbb{U}$.
This implies that
\begin{equation}
|w(z)| = \left|\frac{G(z)-1}{G(z)+1}\right| < 1
\qquad (z\in\mathbb{U}). \label{m1thm2eq4}
\end{equation}

It follows from (\ref{m1thm2eq4}) that
$$
\Re(G(z)) = \Re\left(\frac{zf'(z)}{f(z)}\right)>0 \qquad (z\in\mathbb{U})
$$
so that,
$f(z) \in \mathcal{S^{*}}$.

Spending the same manner with (\ref{m1thm2eq1}),
we conclude $f(z) \in \mathcal{S^{*}}$ for the condition (\ref{m1thm2eq2}).
\end{proof}

\

\begin{thm} \label{m1thm3} \quad
If $f(z)\in\mathcal{A}$ satisfies
$$
\left|\frac{zf'(z)}{f(z)}-1\right|^{\beta} \left|z\left(\frac{zf'(z)}{f(z)}\right)'\right|^{\gamma}<\left(\dfrac{1}{2}\right)^{\gamma}(1-\alpha)^{\beta+\gamma}
\qquad (z\in\mathbb{U})
$$
for some real $0 \leqq \alpha <1$, $\beta$ and $\gamma$ with $\beta + 2\gamma \geqq 0$,
then $f(z) \in \mathcal{S^{*}}(\alpha)$.
\end{thm}

\

\begin{proof}\quad
Defining the function $w(z)$ in $\mathbb{U}$ by
$$
G(z)=\frac{zf'(z)}{f(z)}=\frac{1+(1-2\alpha)w(z)}{1-w(z)}
\qquad (w(z) \neq 1),
$$
we have that $w(z)$ is analytic in $\mathbb{U}$ and $w(0)=0$.
Since
$$
\left(\frac{zf'(z)}{f(z)}\right)'=\frac{2(1-\alpha)w'(z)}{(1-w(z))^{2}},
$$
we obtain that
\begin{align*}
\left|\frac{zf'(z)}{f(z)}-1\right|^{\beta} \left|z\left(\frac{zf'(z)}{f(z)}\right)'\right|^{\gamma}
&= \left|\frac{2(1-\alpha)w(z)}{1-w(z)}\right|^{\beta} \left|\frac{2(1-\alpha)zw'(z)}{(1-w(z))^{2}}\right|^{\gamma}\\
&< \left(\dfrac{1}{2}\right)^{\gamma}(1-\alpha)^{\beta+\gamma}
\qquad (z\in\mathbb{U}).
\end{align*}

If there exists a point $z_{0} \in \mathbb{U}$ such that
$$
\max_{\left| z \right| \leqq \left| z_{0} \right|} \left| w(z) \right|
= \left| w(z_{0}) \right|
= 1,
$$
then Lemma \ref{jack} gives us that $w(z_{0})=e^{i \theta}$ and $z_{0}w'(z_{0})=kw(z_{0})$ $(k \geqq 1)$.

Thus we have
\begin{align*}
\left|\frac{z_{0}f'(z_{0})}{f(z_{0})}-1\right|^{\beta} \left|z_{0}\left(\frac{z_{0}f'(z_{0})}{f(z_{0})}\right)'\right|^{\gamma}
&= \left|\frac{2(1-\alpha)w(z_{0})}{1-w(z_{0})}\right|^{\beta} \left|\frac{2(1-\alpha)z_{0}w'(z_{0})}{(1-w(z_{0}))^{2}}\right|^{\gamma}\\
&= \frac{2^{\beta+\gamma}k^{\gamma}}{|1-w(z_{0})|^{\beta+2\gamma}}(1-\alpha)^{\beta+\gamma}\\
&\geqq \left(\frac{k}{2}\right)^{\gamma}(1-\alpha)^{\beta+\gamma}\\
&\geqq \left(\frac{1}{2}\right)^{\gamma}(1-\alpha)^{\beta+\gamma}.
\end{align*}

This contradicts our condition in the theorem.
Therefore, there is no $z_{0} \in \mathbb{U}$ such that $|w(z_{0})| = 1$.
This means that $|w(z)| < 1$ for all $z \in \mathbb{U}$.
This implies that
\begin{equation}
|w(z)| = \left|\frac{G(z)-1}{G(z)+(1-2\alpha)}\right| <1
\qquad (z \in \mathbb{U}). \label{m1thm3eq1}
\end{equation}

From (\ref{m1thm3eq1}),
we obtain
$$
\Re(G(z))=\Re\left(\frac{zf'(z)}{f(z)}\right)>\alpha
\qquad (z \in \mathbb{U}),
$$
so that, $f(z) \in \mathcal{S^{*}}(\alpha)$.
\end{proof}

\

\section{Conditions for the class \texorpdfstring{$\mathcal{STS}(\mu)$}{STS(mu)}}

\

Using Lemma \ref{jack},
we show the following result for the class $\mathcal{STS}(\mu)$.

\

\begin{thm} \label{m1thm4} \quad
If $f(z)\in\mathcal{A}$ satisfies
$$
\left|\frac{zf'(z)}{f(z)}\right|^{\alpha} \left|z\left(\frac{zf'(z)}{f(z)}\right)'\right|^{\beta} < \left(\dfrac{1}{2}\mu\right)^{\beta}
\qquad (z\in\mathbb{U})
$$
for some real $\alpha \geqq 0$, $\beta >0$ and $\mu = \dfrac{\beta}{\alpha+\beta}$,
then $f(z) \in \mathcal{STS}(\mu)$.
\end{thm}

\

\begin{proof}\quad
Letting
$$
G(z)=\frac{zf'(z)}{f(z)}=\left(\frac{1+w(z)}{1-w(z)}\right)^{\mu}
\qquad (w(z)\neq1)
$$
with $\mu=\dfrac{\beta}{\alpha+\beta}$, we see that $w(z)$ is analytic in $\mathbb{U}$ and $w(0)=0$.
Noting that
$$
\left(\frac{zf'(z)}{f(z)}\right)' = \frac{2\mu w'(z)}{(1-w(z))^{2}}\left(\frac{1+w(z)}{1-w(z)}\right)^{\mu-1},
$$
we have
\begin{align*}
\left|\frac{zf'(z)}{f(z)}\right|^{\alpha} \left|z\left(\frac{zf'(z)}{f(z)}\right)'\right|^{\beta}
&= \left|\frac{1+w(z)}{1-w(z)}\right|^{\alpha\beta+\beta(\mu-1)} \left|\frac{2\mu zw'(z)}{(1-w(z))^{2}}\right|^{\beta}\\
&= \left|\frac{2\mu zw'(z)}{(1-w(z))^{2}}\right|^{\beta} < \left(\dfrac{1}{2}\mu\right)^{\beta}
\qquad (z \in \mathbb{U}).
\end{align*}

If there exists a point $z_{0} \in \mathbb{U}$ such that
$$
\max_{\left| z \right| \leqq \left| z_{0} \right|} \left| w(z) \right|
= \left| w(z_{0}) \right|
= 1,
$$
then Lemma \ref{jack} gives us that $w(z_{0})=e^{i \theta}$ and $z_{0}w'(z_{0})=kw(z_{0})$ $(k \geqq 1)$.

Thus we have
\begin{align*}
\left|\frac{z_{0}f'(z_{0})}{f(z_{0})}\right|^{\alpha} \left|z_{0}\left(\frac{z_{0}f'(z_{0})}{f(z_{0})}\right)'\right|^{\beta}
&= \left|\frac{2\mu z_{0}w'(z_{0})}{(1-w(z_{0}))^{2}}\right|^{\beta}\\
&= \frac{2^{\beta}k^{\beta}\mu^{\beta}}{|1-w(z_{0})|^{2\beta}}\\
&\geqq \left(\frac{k}{2}\mu\right)^{\beta} \geqq \left(\frac{1}{2}\mu\right)^{\beta},
\end{align*}
which contradicts our condition in the theorem.
Therefore, there is no $z_{0} \in \mathbb{U}$ such that $\left|w(z_{0})\right| = 1$.
This means that $\left|w(z)\right| < 1$ for all $z \in \mathbb{U}$.
It follows that
\begin{equation}
G(z)=\left(\frac{1+w(z)}{1-w(z)}\right)^{\mu}
\qquad (|w(z)|<1). \label{m1thm4eq1}
\end{equation}

From (\ref{m1thm4eq1}),
we obtain that
$$
\Re\left(G(z)^{\frac{1}{\mu}}\right) = \Re\left(\left(\frac{zf'(z)}{f(z)}\right)^{\frac{1}{\mu}}\right) > 0
\qquad (z \in \mathbb{U}),
$$
that is, that $f(z) \in \mathcal{STS}(\mu)$.
\end{proof}

\

\end{document}